\documentclass[12pt]{amsart}
\usepackage{amsfonts}
\usepackage{latexsym}
\usepackage{amssymb}
\usepackage{amsmath}
\renewcommand{\H}{\mathbb{H}}

\newcommand{\R}{\mathbb{R}}

\newcommand{\cH}{\mathcal{H}}

\newcommand{\cL}{\mathcal{L}}

\newcommand{\cJ}{\mathcal{J}}

\newcommand{\cS}{\mathcal{S}}

\newcommand{\ph}{\varphi}

\newcommand{\sm}{\setminus}

\newcommand{\res}{\mbox{\LARGE{$\llcorner$}}}
\newcommand{\lan}{\langle}
\newcommand{\ran}{\rangle}

\newcommand{\lra}{\longrightarrow}

\newcommand{\der}{\partial}
\newcommand{\Htwe}{{\mathcal H}^2_{|\cdot|}}

\newtheorem{The}{Theorem}%[section]
\newtheorem{Lem}{Lemma}
\newtheorem{Rem}{Remark}
\newtheorem{Def}{Definition}

\newtheorem{Cor}{Corollary}

\begin{document}

\title
[Nonexistence of horizontal Sobolev surfaces in the Heisenberg group]
{{\bf Nonexistence of horizontal Sobolev surfaces in the Heisenberg group}}
\author{Valentino Magnani}
\address{Valentino Magnani, Dipartimento di Matematica \\
Largo Bruno Pontecorvo 5 \\ 56127, Pisa, Italy}
\email{magnani@dm.unipi.it}
\begin{abstract}
Involutivity is a well known necessary condition for integrability
of smooth tangent distributions. We show that this condition 
is still necessary for integrability with Sobolev surfaces.
We specialize our study to the left invariant
horizontal distribution of the first Heisenberg group $\H^1$. Here we answer a question 
raised in a paper by Z.M.Balogh, R.Hoefer-Isenegger, J.T.Tyson.
\end{abstract}
\maketitle

The Heisenberg group $\H^1$ can be represented as $\R^3$, equipped with
the couple of left inviariant vector fields
\[
X_1(x)=\der_{x_1}-x_2\der_{x_3}\qquad X_2(x)=\der_{x_2}+x_1\der_{x_3}
\]
with respect to the group operation $x\, y=x+y+(0,0,x_1y_2-x_2y_1)$,
for every $x,y\in\R^3$.
In the sequel, we will use the standard Euclidean norm $|\cdot|$ on $\H^1$,
especially when we consider 2-rectifiable sets, in the Federer sense.
We denote by $\cH^\alpha_{|\cdot|}$ the $\alpha$-dimensional Hausdorff 
measure with respect to $|\cdot|$. 
Sobolev mappings with values in $\H^1$ are thought of as having values in 
$\R^3$.
The measures $\cS^\alpha$ and $\cH^\alpha$ are assumed to be contructed 
with respect to a fixed left invariant homogeneous distance of the 
Heisenberg group.
\begin{Def}\label{Sobsur}{\rm We say that a countably $\cH^2_{|\cdot|}$-rectifiable set 
$S$ of $\H^1$ is a {\em Sobolev surface} if it can be written,
up to $\cH^2_{|\cdot|}$-negligible sets, 
as the countable union of graphs of precisely represented Sobolev 
functions of class $W^{1,1}_{loc}$ and defined on open sets of $\R^2$.
}\end{Def}
\begin{Rem}{\rm In view of recent results by J. Mal\'y, D. Swanson and 
W. P. Ziemer, graphs of precisely represented functions in our 
assumptions are 2-rectifiable, see \cite{MSZ03}.
Then the hypothesis of rectifiability in Definition~\ref{Sobsur} could be removed.
}\end{Rem}
The distribution of admissible directions in the Heisenberg group
is given by the following {\em horizontal subspaces}
\[
H_y\H^1=\{\lambda_1 X_1(y)+\lambda_2X_2(y)\mid \lambda_j\in\R\}
\quad\mbox{for every} \quad y\in\H^1.
\]
The collection of all horizontal subspaces $H_y\H^1$, $y\in\H^1$,
seen as subbundle of $T\H^1$ is the so-called {\em horizontal subbundle}
and it is denoted by $H\H^1$.
\begin{Def}{\rm A 2-rectifiable set $S$ in $\H^1$ is {\em horizontal} if
for $\cH^2_{|\cdot|}$-a.e. $y\in S$ we have $\mbox{Tan}(S,y)\subset H_y\H^1$.
We also say that $S$ is $\Htwe$-a.e. tangent to $H\H^1$.
}\end{Def}
\begin{Rem}{\rm Then nonexistence of horizontal smooth 2-dimensional submanifolds
in $\H^1$ is a simple consequence of the fact that the horizontal distribution 
given by horizontal subspaces $H_y\H^1$ is non-involutive.
In fact, $[X_1,X_2]=2\der_{y_3}$ and this vector field clearly is not a 
linear combination of $X_1$ and $X_2$.
}\end{Rem}
\begin{Lem}\label{nnex1}
Let $f\in W_{loc}^{1,1}(\Omega,\H^1)$ be a graph parametrization of
a Sobolev surface. Then the following system
\begin{eqnarray}\label{cnteqsHeis}
\left\{\begin{array}{l}
f^3_{x_1}=f^1f^2_{x_1}-f^2f^1_{x_1}\\
f^3_{x_2}=f^1f^2_{x_2}-f^2f^1_{x_2}
\end{array}\right.
\end{eqnarray}
fails to hold in a subset of positive measure.
\end{Lem}
{\sc Proof.}
Recall that $\Omega$ is an open subset of $\R^2$.
We can rewrite the system \eqref{cnteqsHeis} in terms of 
differential forms as the a.e. pointwise validity of 
\[
df^3=f^1df^2-f^2df^1\,.
\]
Since $f$ parametrizes a graph, it can be represented in three possible ways,
where it always happens that either $f^1$ or $f^2$ is a coordinate function.
Thus, one of these components clearly is in $W^{1,1}_{loc}(\Omega)$
and the remaining one is smooth.
As a consequence, both $f^1df^2$ and $f^2df^1$
can be weakly differentiated and the weak exterior differential satisfies the
formula
\[
d\big(f^1df^2-f^2df^1\big)=2\,df^1\wedge df^2\,.
\]
Clearly, $d(df^3)=0$ in the distributional sense, hence 
\[
\int_\Omega *\big(df^1\wedge df^2\big)\;\phi\;d\cL^2=0
\]
for $\phi\in C_c^\infty(\Omega)$, where $*\big(df^1\wedge df^2\big)=\det(f^i_{x_j})$.
We have proved that $\nabla f^1(x)$ and $\nabla f^2(x)$ are not linearly independent
for a.e. $x\in\Omega$. Due to \eqref{cnteqsHeis}, it follows that
the rank of $\nabla f(x)$ is less than or equal to one for a.e. $x\in\Omega$. 
This conflicts with the fact that $f$ parametrizes a graph. $\Box$
\begin{Rem}\label{weakextdif}{\rm
In the previous proof we have used the notion of weak exterior differential
of a locally summable $k$-form $\alpha$ on an open set $\Omega$ of $\R^n$.
Recall that the locally summable $(k+1)$-form $\beta$ is the weak exterior differential of
$\alpha$ if for every smooth compactly supported $(n-k-1)$-form $\phi$, we have
\[
\int_\Omega\lan\alpha,*\,d\phi\ran\,d\cL^n
=(-1)^{k+1}\int_\Omega\lan\beta,\,*\phi\ran\,d\cL^n
\]
Here $*$ denotes the Hodge operator with respect to the volume
form $dx_1\wedge\cdots\wedge dx_n$. Notice also that $\beta$ is uniquely defined.
The validity of formulae $d(f^1\,df^2)=df^1\wedge df^2$ and
$d(f^2\,df^1)=df^2\wedge df^1$ used in the previous proof can be obtained
by standard smooth approximation arguments.
}\end{Rem}
\begin{Rem}\label{htan}{\rm
One can check that the pointwise validity of
\eqref{cnteqsHeis} coincides with the pointwise validity of 
either $df(x)(T_x\R^2)\subset H_{f(x)}\H^1$ or equivalently 
Tan$(S,f(x))=H_{f(x)}\H^1$, where $S$ is parametrized by $f$. 
}\end{Rem}
\begin{The}
There do not exist horizontal Sobolev surfaces in $\H^1$.
\end{The}
{\sc Proof.}
By contradiction, we assume that $\Sigma$ is a horizontal Sobolev surface in $\H^1$.
Then we have $f\in W^{1,1}_{loc}(\Omega,\H^1)$ that is the graph of some
$W^{1,1}_{loc}$-function and such that $f(\Omega)$ is $\Htwe$-a.e. tangent to $H\H^1$.
If we could find a set $E\subset\Omega$ of positive measure where \eqref{cnteqsHeis} 
fails to hold, then by Theorem~1.2 of \cite{MSZ03} and in view of Remark~\ref{htan},
we would get a subet $f(E)\subset\Sigma$ of positive measure $\Htwe$
that is a.e. not tangent to $H\H^1$. This conflicts with our assumption
of horizontality, hence we have proved that \eqref{cnteqsHeis} holds a.e. in $\Omega$.
The latter assertion conflicts with Lemma~\ref{nnex1} and concludes the proof. $\Box$
\begin{Rem}{\rm
Notice that, by definition, each Sobolev surface has positive 
measure $\Htwe$, hence one immediately observes that it also
has positive measure $\cH^2$. On the other hand, $\Htwe$-negligible sets
cannot have positive measure $\cH^3$, since this
measure is absolutely continuous with respect to $\Htwe$, as it has
been shown in \cite{BRSC03}.
}\end{Rem}
\begin{The}\label{SobH3}
Every Sobolev surface $\Sigma\subset\H^1$ satisfies $\cH^3(\Sigma)>0$. 
\end{The}
{\sc Proof.}
By definition of Sobolev surface, we can find a precisely
represented function $u\in W^{1,1}_{loc}(\Omega)$, where $\Omega$ is an open subset
of $\R^2$, such that the graph of $u$ is contained in $\Sigma$. 
Suppose that the graph is of the form 
\begin{eqnarray}\label{graph1}
\Omega\ni(x_1,x_2)\lra f(x_1,x_2)=\big(u(x_1,x_2),x_1,x_2\big).
\end{eqnarray}
Lemma~\ref{nnex1} ensures that there is a subset $A\subset\Omega$ of
positive measure such that
\begin{eqnarray}\label{nothold}
\left\{\begin{array}{l}
u-x_1u_{x_1}=0 \\
1+x_1u_{x_2}=0
\end{array}\right.
\quad\mbox{does not hold at every point of $A$.}
\end{eqnarray}
Taking into account the classical Whitney extension theorem, 
see for instance 3.1.15 of \cite{Fed} and the lemma of Section~3 in
\cite{Haj93}, one can find a 
bounded subset with positive measure $A_0\subset A$ and a $C^1$
smooth function $v:\R^2\lra\R$ such that $u$ is everywhere differentiable in
$A_0$ and there coincides with $u$ along with its gradient.
We define the submanifold
\[
\Sigma_1=\{(y_1,y_2,y_3)\in\H^1\mid v(y_2,y_3)-y_1=0,\,
(y_2,y_3)\in\Omega\}\,.
\]
Taking into account formulae (5.1) and (5.2) of \cite{Bal} for $n=1$, we have
\begin{equation}\label{s3}
\cS^3\res\Sigma_1=|{\bf n}_H|\,d\cH^2_{|\cdot|}\res\Sigma_1,
\end{equation}
where $\cS^3$ is the spherical Hausdorff measure with respect to a
fixed Heisenberg metric. The length of the horizontal normal with respect
to the Euclidean metric is given by
\[
|{\bf n}_H(v(y),y)|^2=(1+y_2v_{y_3}(y))^2+(v_{y_2}(y)+v(y)v_{y_3}(y))^2\,
\]
since it is equal to
$\lan{\bf n}(v(y),y),X_1(n(v(y),y))\ran^2+\lan{\bf n}(n(v,y)),X_2(n(v,y))\ran^2$,
where we have set $y=(y_2,y_3)$.
Taking into account \eqref{nothold}, for every $y\in A_0$, we have that
$v(y)=u(y)$ and either
\[
|v(y)-y_2v_{y_2}(y)|>0\quad\mbox{or}\quad |1+y_2v_{y_3}(y)|>0\,.
\]
If $(1+y_2v_{y_3}(y))\neq0$ on a subset $E\subset A_0$ of positive measure,
then $|{\bf n}_H(v(y),y)|>0$ for every $y\in E$.
By Theorem~1.2 of \cite{MSZ03}, $f$ preserves $\Htwe$-negligible sets
and also $\Htwe(f(E))>0$.
As a result, due to \eqref{s3} we get $\cS^3(f(E))>0$, where $f(E)\subset\Sigma$.
The remaining case is that $1+y_2v_{y_3}(y)=0$ for a.e. $y\in A_0$.
In particular, $y_2\neq0$ and $|v(y)-y_2v_{y_2}(y)|>0$ for a.e. $y\in A_0$.
As a consequence, 
\[
0<|y_2v_{y_2}(y)-v(y)|=|y_2|\,|v_{y_2}(y)+v(y)v_{y_3}(y)|
\leq|y|\,|{\bf n}_H(v(y),y)|
\]
for a.e. $y\in A_0$. Thus, arguing as before, we get $\cS^3(f(A_0))>0$,
where $f(A_0)\subset\Sigma$.
This concludes the proof in the case the graph has the form \eqref{graph1}.
The remaining two cases have analogous proof. $\Box$
\begin{Cor}
There do not exist Sobolev surfaces $\Sigma$ in $\H^1$ such that
$0<\cH^2(\Sigma)<\infty$.
\end{Cor}
This corollary answers a question raised in \cite{BHIT} by
Z. M. Balogh, R. Hoefer-Isenegger and J. T. Tyson about the possibility to
construct sets with finite and positive measure $\cH^2$ with regularity
between BV and Lipschitz. The authors show that 
there exist graphs of BV functions that have this property, although 
this is not true for Lipschitz parametrizations, as it has been shown
in \cite{AmbKir} by L. Ambrosio and B. Kirchheim.
Precisely, Lipschitz parametrizations from $\R^2$ to $\H^1$
are considered with respect to the Carnot-Carath\`eodory distance of $\H^1$
and this also implies the local Lipschitz property with respect
to the Euclidean distance fixed in $\H^1$. Here we wish to mention that
Lipschitz maps between stratified groups a.e. satisfy their associated contact equations,
\cite{Mag10}, and these equations in our case exactly correspond to the system \eqref{cnteqsHeis}.
\begin{Rem}{\rm
Notice that the previous lemma precisely shows that the closure of
the set where \eqref{cnteqsHeis} fails to hold coincides with $\Omega$.
On the other hand, it is still possible to construct even $C^{1,\alpha}$
parametrizations of graphs in $\H^1$, with $0<\alpha<1$,
where \eqref{cnteqsHeis} holds in a subset of positive measure, \cite{Bal}. 
Then this subset must have empty interior.
}\end{Rem}
It is natural to consider our previous results for parametrized surfaces,
that are not necessarily graphs.
In fact, one can extend the notion of
Sobolev surface to suitable images of Sobolev mappings. 
Clearly, this is a weaker notion than the previous one.
\begin{Def}\label{Sobpar}{\rm
We say that a countably $\cH^2_{|\cdot|}$-rectifiable set 
$S$ of $\H^1$ is a {\em parametrized $W^{1,p}$-Sobolev surface},
if it can be written, up to $\cH^2_{|\cdot|}$-negligible sets, 
as the countable union of images of $W^{1,p}_{loc}$-Sobolev mappings
on open subsets of $\R^2$, that sends $\cH^2_{|\cdot|}$-negligible sets 
into $\cH^2_{|\cdot|}$-negligible sets and that have a.e. maximal rank.
}\end{Def}
Although in the previous definition rectifiability is a consequence
of the assumptions on the Sobolev parametrizations, we have preferred
to stress this important property.
\begin{Rem}{\rm
In Definition~\ref{Sobpar}, we have assumed also a sort of Lusin's condition on
the parametrization, namely, that of preserving $\cH^2_{|\cdot|}$-negligible sets.
This is an important assumption, since one can find for instance
Sobolev mappings of $W^{1,2}(\R^2,\R^3)$
whose image coincides with all of $\R^3$,
see \cite{HajTys} for more general results in this vein.
}\end{Rem}
\begin{Rem}{\rm
It is also natural to assume that the Sobolev parametrizations 
considered in Definition~\ref{Sobpar} have a.e. maximal rank. In fact, without
this assumption one can consider the smooth mapping
$\{x\in\R^2\mid 0<|x|<1\}\ni x\lra(0,0,|x|)\in\H^1$ 
whose image has positive and finite measure $\cH^2$, as was
already pointed out in \cite{BHIT}.
}\end{Rem}
\begin{Lem}\label{nnex2}
Let $f\in W_{loc}^{1,4/3}(\Omega,\H^1)$ be a Sobolev mapping with
a.e. maximal rank. Then conclusions of Lemma~\ref{nnex1} still hold.
\end{Lem}
{\sc Proof.}
Suppose by contradiction that \eqref{cnteqsHeis} holds a.e. in $\Omega$,
then it can be written as follows
\[
df^3=f^1df^2-f^2df^1\,.
\]
Then the weak exterior differential $f^1df^2-f^2df^1$ is clearly vanishig
and equals twice the distributional Jacobian, see Section~7.1 of \cite{IM}.
Thus, we have
\[
\lan\cJ_F,\ph\ran=-\int_\Omega f^1\,d\ph\wedge df^2=0\quad
\mbox{for every $\ph\in C_c^\infty(\Omega)$}\,,
\]
where we have set $F=(f^1,f^2):\Omega\lra\R^2$.
Then we apply Lemma~7.1.1 of \cite{IM} to get that for a.e. $x\in\Omega$
there exists the limit
\[
\lim_{t\to0^+}\cJ_j*\Phi_t(x)=J(x,F)\,.
\]
Where $J(x,F)=\det\big((f^i_{x_j})_{i,j=1,2}\big)$ is the pointwise Jacobian.
Since the distributional Jacobian is vanishing, we have that $J(x,F)=0$
for a.e. $x\in\Omega$. Taking into account \eqref{cnteqsHeis}, 
we have proved that the rank of $f$ is a.e. less than or equal to one.
This conflicts with our assumptions on $f$. $\Box$
\begin{Rem}{\rm
The previous lemma relies on the notion of distributional jacobian
and its properties. We address the reader to the recent monograph
\cite{IM} for a thorough presentation of this topic along with a number
of related arguments.}
\end{Rem}
\begin{Lem}
Let $f\in W^{1,1}_{loc}(\Omega,\R^n)$, where $\Omega\subset\R^k$ is an open subset
and $k\leq n$ and suppose that $f$ preserves $\cH^k_{|\cdot|}$-negligible sets.
Then the following area formula holds
\begin{equation}\label{areaf}
\int_E Jf(x)\,dx=\int_{\R^n}\,N_f(y,E)\,d\cH^k_{|\cdot|}(y)\,,
\end{equation} 
where $E$ is a measurable set in $\Omega$ and $Jf(x)$ denotes the 
jacobian of the approximate differential of $f$ at $x$.
\end{Lem}
{\sc Proof.}
One argues as in \cite{Haj93}. In fact, the area formula
holds for Lipschitz mappings and it is possible to find an 
increasing sequence of measurable sets $X_k$,
whose union gives $E$ up to an $\cH^k_{|\cdot|}$-negligible set
and such that $f_{|X_k}$ is Lipschitz.
By our assumption $\cH^k_{|\cdot|}\big(f(E\sm\cup_k X_k)\big)=0$,
hence Beppo-Levi convergence theorem concludes the proof. $\Box$
\begin{The}
There do not exist horizontal parametrized $W^{1,p}$-Sobolev surfaces in $\H^1$
for every $p\geq 4/3$.
\end{The}
{\sc Proof.} 
Let $\Sigma$ be a parametrized $W^{1,4/3}$-Sobolev surface.
By definition we can find a Sobolev mapping $f\in W^{1,4/3}(\Omega,\H^1)$
with a.e. maximal rank that sends $\cH^2_{|\cdot|}$-negligible sets
into $\cH^2_{|\cdot|}$-negligible sets and such that $f(\Omega)\subset\Sigma$.
By contradition, suppose that $\Sigma$ is horizontal.
We wish to prove that $df(x)(T_x\R^2)\subset H_{f(x)}\H^1$ 
a.e. in $\Omega$. In fact, if this were not the case, then one could find a set
of positive measure $E$ in $\Omega$ where the previous condition of
horizontality does not hold. By area formula \eqref{areaf} and the hypothesis
on the rank of $f$ we would get $\Htwe(f(E))>0$, where 
$f(E)$ is not tangent to $H\H^1$ at $\Htwe$-a.e. point.
This conflicts with our hypothesis on $\Sigma$.
Since $df(x)(T_x\R^2)\subset H_{f(x)}\H^1$ is equivalent to the validity
of \eqref{cnteqsHeis} at $x$, we have proved the a.e. validity
of \eqref{cnteqsHeis} in $\Omega$ and this conflicts with Lemma~\ref{nnex2}. $\Box$ 
\begin{The}
Let $p\geq 4/3$ and let $\Sigma$ be a parametrized $W^{1,p}$-Sobolev surface in $\H^1$.
Then $\cH^3(\Sigma)>0$. 
\end{The}
{\sc Proof.}
By hypothesis, we have a $W^{1,p}_{loc}$-mapping $f:\Omega\lra\H^1$ with
a.e. maximal rank that sends $\Htwe$-negligible sets into $\Htwe$-negligible sets
and such that $f(\Omega)\subset\Sigma$.
By Lemma~\ref{nnex2}, one can find a set $E\subset\Omega$ of positive measure
such that $f$ is everywhere approximately differentiable on $E$ and 
the system \eqref{cnteqsHeis} fails to hold everywhere on this set.
It is also not restrictive assuming that each point of $E$ is a density point
and the approximate differential has maximal rank.
Again, by Whitney extension theorem and the lemma of Section~3 in
\cite{Haj93}, one can find a subset $E_0$ of $E$ with positive measure and 
a $C^1$ mapping $g:\R^2\lra\H^1$ such that $g_{|E_0}=f_{|E_0}$ and
the approximate differential of $f$ along with the differential of $g$ coincide on $E_0$.
Let us fix $x_0\in E_0$ and notice that for a fixed $r_0>0$ sufficiently small
we have both $\cL^2(B_{x_0,r_0}\cap E_0)>0$
and $g(B_{x_0,r_0})=\Sigma_0\subset\H^1$ is an embedded surface.
Up to possibly shrinking $S_0$, it is not restrictive to assume that
it is a graph around $x_0$. As in the proof of Theorem~\ref{SobH3}, we apply \eqref{s3},
then getting
\begin{equation}\label{s3a}
\cS^3\res\Sigma_0=|{\bf n}_H|\,d\cH^2_{|\cdot|}\res\Sigma_0.
\end{equation}
By properties of $g$ and \eqref{areaf}, 
we have $S_0=f(B_{x_0,r_0}\cap E_0)\subset\Sigma_0\cap\Sigma$
and $\Htwe(S_0)>0$.
Since \eqref{cnteqsHeis} does not hold on $E_0$, then 
${\bf n}_H(f(x))\neq0$ for every $x\in E_0$, hence
\eqref{s3a} gives $\cS^3(S_0)>0$. This concludes the proof. $\Box$
\begin{Cor}
For every $p\geq 4/3$, there do not exist parametrized $W^{1,p}$-Sobolev surfaces 
$\Sigma$ such that $0<\cH^2(\Sigma)<\infty$.
\end{Cor}
\vskip.25cm
As a final comment, we wish to point out how this note leaves open
the question about existence of horizontal parametrized
$W^{1,p}$-Sobolev surfaces with $1\leq p<4/3$.
\vskip.25cm
{\bf Acknowledgements.} It is a great pleasure to thank Tadeusz Iwaniec
for pleasant discussions and for his kind suggestion about Lemma~\ref{nnex2}.

\end{document}